\documentclass{article}

\usepackage{amsmath}
\usepackage{amssymb}
\usepackage{latexsym}
\usepackage{amsxtra}
\usepackage{amscd}
\usepackage{theorem}
\usepackage{xypic}

\theoremstyle{plain}

{\theorembodyfont{\slshape}

        \newtheorem{thm}{Theorem}[section]
        \newtheorem{cor}[thm]{Corollary}
        \newtheorem{lem}[thm]{Lemma}
        \newtheorem{prop}[thm]{Proposition}

}

{\theorembodyfont{\rmfamily}

        \newtheorem{defn}[thm]{Definition}
        
        \newtheorem{rem}[thm]{Remark}
        \newtheorem{exa}[thm]{Example}

}

\renewcommand{\em}{\sl}

\newcommand{\proof}{{\bf Proof:\ }}
\newcommand{\Endproof}{\hspace*{\fill} $\Box$ \vspace{1ex} \noindent }

\makeatletter
\renewcommand{\subsubsection}{\@startsection{subsubsection}{3}%
        {\z@}{-3.25ex plus -1ex minus-.2ex}{-1em}{\bf}}
\makeatother

\newcommand{\PP}{\mathbb{P}}
\newcommand{\ZZ}{\mathbb{Z}}

\newcommand{\RR}{\mathbb{R}}
\newcommand{\QQ}{\mathbb{Q}}

\renewcommand{\AA}{\mathbb{A}}

\newcommand{\F}{\mathcal{F}}

\newcommand{\W}{\mathcal{W}}
\newcommand{\U}{\mathcal{U}}

\newcommand{\X}{\mathcal{X}}
\newcommand{\Y}{\mathcal{Y}}

\newcommand{\J}{\mathcal{J}}

\newcommand{\Xs}{\mathsf{X}}
\newcommand{\Ys}{\mathsf{Y}}

\newcommand{\OO}{\mathcal{O}}

\newcommand{\ac}{^{\rm ac}}

\newcommand{\cusp}{^{\rm csp}}

\newcommand{\Kh}{\hat{K}}

\newcommand{\QQb}{\bar{\QQ}}

\newcommand{\Hom}{{\rm Hom}}

\newcommand{\Aut}{{\rm Aut}}

\newcommand{\Ker}{{\rm Ker}}

\newcommand{\GL}{{\rm GL}}

\newcommand{\Spec}{{\rm Spec\,}}
\newcommand{\Spf}{{\rm Spf}}

\newcommand{\red}{{\rm red}}

\newcommand{\inj}{\hookrightarrow}
\newcommand{\To}{\;\longrightarrow\;}
\newcommand{\iso}{\stackrel{\sim}{\to}}

\newcommand{\pfeil}[1]{\stackrel{#1}{\to}}

\newcommand{\abs}[1]{\lvert#1\rvert}

\makeatletter
\renewcommand{\subsection}{\@startsection{subsection}{2}%
        {\z@}{-3.25ex plus -1ex minus-.2ex}{-1em}{\bf}}
\makeatother

\begin{document}

\title{Some remarks on open analytic curves over non-archimedian fields}

\author{Stefan Wewers}

\date{}

\maketitle

\begin{abstract}
  We study open analytic curves over non-archimedian fields and their
  formal models. In particular, we give a criterion, in terms of
  \'etale cohomology, when such a formal model is (almost)
  semistable.\footnote{Part of this work has been achieved thanks to
    the support of the European Commission through its 6th Framework
    Program "Structuring the European Research Area" and the
    contract Nr. RITA-CT-2004-505493 for the provision of
    Transnational Access implemented as Specific Support Action.}
\end{abstract}

\section*{Introduction}

Let $\Xs$ be a smooth analytic curve over a non-archimedian complete
valued field $K$. A {\em formal model} of $\Xs$ is a formal scheme
$\X$ over the valuation ring of $K$ with generic fiber $\Xs$. The
purpose of this note is to formulate a criterion which is, in certain
concrete situations, able to decide whether the model $\X$ is
semistable.

The chief motivation to formulate such a criterion comes from the
author's work on the semistable reduction of Lubin-Tate spaces of
dimension one \cite{LT}. There the analytic curve $\Xs$ in question
arises as a finite \'etale Galois cover of the open unit disk,
\[
      f:\Xs\to\Ys := \{\,y\in K \,\mid\, \abs{y} < 1 \,\}.
\]
Such a cover belongs to a class of analytic spaces which we call {\em
  open analytic curves} and which provides a non-archimedian analogue
of open Riemann surfaces with finitely many holes. Coleman has studied
this class of analytic spaces (which he calls {\em wide open spaces}),
also in connection with the problem of semistable reduction
(\cite{Coleman03}). One difference to our approach is that in
\cite{Coleman03} the space $\Xs$ is always considered as an open
analytic subspace of an algebraic curve $C$ (whose semistable
reduction one wants to analyze). For the applications we have in mind
it is however important to consider $\Xs$ as an object of its own,
independent of any embedding into an algebraic curve. This does not
result in any serious problems but, since open analytic curves are not
quasi-compact, one has to be somewhat careful in applying results
which rely on finiteness arguments. For instance, it has seemed safer
to the author to always work over a {\em discrete} valued field $K$.
Moreover, lacking adequate references it seemed necessary to
reformulate and discuss many definitions which are well known in the
algebraic context. This is done in Section 1.

In Section 2 we look at the $\ell$-adic cohomology of open analytic
curves and at the vanishing cycles on the special fiber of formal
models. Using the framework provided by the work of Berkovich
(\cite{Berkovich93}, \cite{Berkovich96}), we formulate a first version
of our semistable reduction criterion. Given a formal model $\X$ of an
open analytic curve $\Xs$ with special fiber $\X_s$, we say that $\X$
is {\em almost semistable} if for every closed point $z\in\X_s$ the
formal fiber 
\[
    \Xs_z:=]z[_{\X}\subset\Xs
\] 
is an open analytic curve of genus zero, i.e.\ is isomorphic to the
complement of finitely many closed disks lying inside an open disk.
(Note that $\X$ is semistable if and only if $\Xs_z$ is either an open
disk or an open annulus, for all $z$.) Then we prove that $\X$ is
almost semistable if the image of the natural map
\[
      H^1(\X_s) \to H^1(\Xs)
\]
is equal to the {\em cuspidal part} of $H^1(\Xs)$, i.e.\ the image of
cohomology with compact support.

One thing that should be mentioned is that the conclusion of our
criterion is somewhat stronger than stated above. In addition to the
assertion that $\X$ is almost semistable one concludes that $\Xs$
has {\em tree-like reduction}. By this we mean that the special fiber
of any semistable model of $\Xs$ has a graph of components which is a
tree. So our criterion would not work in a situation where this
additional conclusion does not hold. Fortunately, the Lubin-Tate
spaces studied in \cite{LT} do have tree-like reduction, and our
method can be applied to them. 

In Section 3 we give an `equivariant version' of the above criterion.
Let $f:\Xs\to\Ys$ be a finite \'etale Galois cover of open analytic
curves, with Galois group $G$. Let $\Y$ be a semistable model of $\Ys$
and let $\X$ be the normalization of $\Y$ in $\Xs$. The question is
now: is $\X$ (almost) semistable?  Obviously, it suffices to verify
the above criterion on the $\tau$-isotypical part of the cohomology of
$\Xs$, for every irreducible representation $\tau$ of $G$. The
$\tau$-isotypical part of the cohomology of $\Xs$ can be expressed in
terms of the cohomology of the \'etale sheaf $\F_\tau$ on $\Ys$
associated to the Galois cover $f:\Xs\to\Ys$ and the representation
$\tau$. In this situation we prove (Proposition \ref{criterionprop}):

\begin{prop}
  Suppose that for every irreducible $G$-representation $\tau$ there
  exists an affinoid subdomain $U\subset\Ys$ with the following properties:
  \begin{enumerate}
  \item
    $U=]W[_{\Y}$ for an open subset $W$ of the special fiber of $\Y$, and
  \item
    the sheaf $\F_\tau$ is `resolved' over $U$ (see Definition \ref{resolvedef}). 
  \end{enumerate}
  Then the formal model $\X$ is almost semistable.
\end{prop}
  
This equivariant version of our criterion, although it is essentially
a reformulation of the original version, turns out to be very useful.
A psychological advantage is that one `can forget about the curve
$\Xs$'. In practice, one has to find for every representation $\tau$
an affinoid $U$ such that Condition (ii) holds. To do this there are
several useful tools available, for instance Huber's theory of {\em
  Swan conductors} \cite{Huber01}. The $U$'s that arise in this
way determine the semistable model $\Y$ we should take. If $\Ys$ is an
open disk then this last step is trivial.

We end this paper with an extension of the above criterion which
sometimes allows one to conclude that the model $\X$ is actually
semistable (and not just almost semistable). Here we use in an
essential way results and arguments from Raynaud's paper
\cite{Raynaud90}.

\section{Open analytic curves} \label{open}

\subsection{The definition}

We fix, once and for all, a field $K_0$ which is complete with respect
to a discrete non-archimedian valuation $\abs{\,\cdot\,}$, and whose
residue field $k$ is algebraically closed and of positive
characteristic $p>0$. We choose an algebraic closure $K\ac_0$ of $K_0$
and extend the valuation $\abs{\,\cdot\,}$ to $K\ac_0$.

\begin{defn} \label{opendef}
  An {\em open analytic curve} is given by a pair $(K,\Xs)$, where
  $K\subset K\ac_0$ is a finite extension of $K_0$ and $\Xs$ is a
  rigid analytic space over $K$. We demand that $\Xs$ is isomorphic to
  $C-D$, where $C$ is the analytification of a smooth projective curve
  over $K$ and $D\subset C$ is an affinoid subdomain intersecting
  every connected component of $C$.

  A morphism between two open analytic curves $(K_1,\Xs_1)$ and $(K_2,\Xs_2)$
  is an element of the direct limit 
  \[
        \Hom(\Xs_1,\Xs_2):= \varinjlim_{K_3} \Hom(\Xs_1\otimes K_3,\Xs_2\otimes K_3),
  \]
  where $K_3\subset K\ac_0$ ranges over all common finite extensions of $K_1$ and $K_2$.
\end{defn}

Definition \ref{opendef} corresponds to the definition of {\em wide
  opens} in \cite{Coleman03}. However, our definition is more
complicated because we insist on having a field of definition with a
{\em discrete} valuation, whereas Coleman works over $\Kh\ac_0$, the
completion of $K_0\ac$.  Note that a morphism
$\Xs_1\otimes\Kh\ac_0\to\Xs_2\otimes\Kh\ac_0$ of rigid analytic spaces
over $\Kh\ac_0$ may not descend to an element of $\Hom(\Xs_1,\Xs_2)$,
so this difference is more than just formal.

Most facts about open analytic curves that we are interested in only
hold after replacing its field of definition by some finite extension
(e.g.\ the existence of a semistable model).  However, the exact
choice of this extension will not be important for us. Therefore we
will simply write $\Xs$ instead of $(K,\Xs)$ to denote an open
analytic curve. Whenever it is necessary to mention the field $K$
(which we then call the {\em field of definition} of $\Xs$) we will
always assume that it is chosen `sufficiently large'. For instance, if
we say that $\Xs$ is connected we actually mean that $\Xs\otimes K'$
is connected for every finite extension $K'/K$.

\begin{exa} \label{openexa}
\begin{enumerate}
\item
  An {\em open disk} is an open analytic curve isomorphic to the open unit disk
  \[
        D(0,1) = \{\, x \,\mid\, \abs{x}<1 \,\}.
  \]
\item An {\em open annulus}   is an open
  analytic curve isomorphic to the standard annulus
  \[
        A(\epsilon,1) = \{\, x \,\mid\, \epsilon<\abs{x}<1 \,\},
  \]
  for some $\epsilon\in \abs{K^\times}$ with $\epsilon<1$. The number
  $\epsilon$ is easily seen to depend only on the isomorphism class of
  $A(\epsilon,1)$.
\end{enumerate}
\end{exa} 

\subsection{Formal and semistable models}

Let $\Xs$ be an open analytic curve, with field of definition $K$. Let
$\OO$ denote the valuation ring of $K$ and $\wp$ the maximal ideal of
$\OO$. If $Y$ is a (formal) scheme over $\OO$ we shall write
$Y_s:=Y\otimes k$ to denote its special fiber and $Y_\eta:=Y\otimes K$
to denote its generic fiber (whenever this makes sense). For a
constructible subset $Z\subset Y_s$ we write $]Z[_Y\subset Y_\eta$ for
the formal fiber of $Z$ (which is an open rigid analytic subspace).
By a {\em curve} we mean a morphism of schemes which is flat and of
finite type and has geometric fibers of pure dimension one.

\begin{defn} \label{modeldef} An {\em algebraic model} of $\Xs$ is a
  triple $(Y,Z,\varphi_Y)$, where $Y$ is a curve over $\OO$, $Z\subset
  Y_s$ is a reduced closed subset and $\varphi:\Xs\iso ]Z[_Y$ is an
  isomorphism of rigid analytic spaces over $K$. The algebraic model
  $(Y,Z,\varphi_Y)$ is called {\em good} if $Y$ is normal and $Y_s$ is
  reduced. It is called {\em minimal} if it is good and if $Z$ is
  purely of dimension zero. It is called {\em semistable} if $Y_s$ is
  a semistable curve and $Z$ is purely of dimension one. 

  Let $(Y,Z,\varphi_Y)$ be a semistable algebraic model and $W$ an
  irreducible component of $Y_s$ contained in $Z$. We call $W$ {\em
    instable} if it is isomorphic to a projective line and intersects
  other irreducible components of $Y_s$ in at most two points. We call
  the algebraic model $(Y,Z,\varphi_Y)$ {\em stable} if there do not
  exist any instable components of $Y_s$ contained in $Z$.

  A {\em formal model} of $\Xs$ is a pair $(\X,\varphi_{\X})$, where
  $\X$ is a formal scheme over $\OO$ and $\varphi_{\X}:\Xs\iso
  \X_\eta$ is an isomorphism of rigid analytic spaces over $K$,
  satisfying the following condition. There exists an algebraic model
  $(Y,Z,\varphi_Y)$ such that $\X$ is isomorphic to the formal
  completion of $Y$ along $Z$, and the isomorphism $\varphi_{\X}$ is
  induced from $\varphi_Y$. The triple $(Y,Z,\varphi_Y)$ is called an
  {\em algebraization} of $(\X,\varphi_{\X})$. A formal model
  $(\X,\varphi_{\X})$ is said to be {\em good} (resp.\ {\em minimal},
  resp.\ {\em semistable}, resp.\ {\em stable}) if it has an
  algebraization which is good (resp.\ minimal, resp.\ semistable,
  resp.\ stable).
\end{defn}

Whenever this is unlikely to cause confusion, we will omit the isomorphisms
$\varphi_Y$ and $\varphi_{\X}$ from the notation.

\begin{rem} \label{formalrem}
\begin{enumerate}
\item Let $\X$ be a formal model of $\Xs$. Then the formal scheme $\X$
  is {\em special} in the sense of \cite{Berkovich96}, \S 1, and
  therefore the generic fiber $\X_\eta$ is well defined. Furthermore,
  if $(Y,Z)$ is an algebraization of $\X$ then $Z$ can be identified
  with the closed subscheme of $\X_s$ corresponding to an ideal of
  definition of the formal scheme $\X$, and is hence independent of
  $Y$. We call $Z$ the {\em reduction} of $\X$. We denote by 
  \[
      \red_{\X}:\Xs\cong\X_\eta \To Z
  \]
  the reduction map. 
\item Note that the formal scheme $\X$ is {\em not} topologically of
  finite type over $\OO$. This corresponds to the
  fact that the rigid analytic space $\Xs$ is not quasi-compact.
  Note also that $Z\neq \X_s$ and that $\OO_{\X}\cdot\wp\subset\OO_{\X}$ is
  {\em not} an ideal of definition for $\X$.
\item The property of being a good (resp.\ minimal, semistable or
  stable) formal model is stable under finite extension of the base
  field $K$. More precisely, if $K'/K$ is a finite extension, $\OO'$
  the valuation ring of $K'$ and if $\X$ is good (resp.\ minimal,
  semistable or stable) then the formal model $\X':=\X\otimes\OO'$ has
  the same property.
\item Let $Y$ be a proper curve over $\OO$ and $Z\subset Y_s$ a reduced
  closed subset which has nontrivial intersection with every connected
  component of $Y_s$. Then the rigid space $\Xs:=]Z[_Y$ is an open
  analytic curve and $(Y,Z)$ is an algebraic model. Indeed, after
  blowing down all irreducible components of $Y_s$ which do not meet $Z$,
  we may assume that the open subset $V:=Y_s-Z$ is affine and hence
  $D:=]V[_Y$ is an affinoid subdomain of $Y_\eta$.
\item Let $\X$ be a semistable model of $\Xs$. Then there exists a
  semistable algebraization $(Y,Z)$ of $\X$ such that $V:=Y_s-Z$ is
  isomorphic to a disjoint union of affine lines. This follows from a
  standard argument using formal patching, see e.g.\
  \cite{FresnelvdPut}. As a consequence we can write $\Xs=C-D$ where
  $C$ is a smooth projective curve over $K$ and $D\subset C$ is an
  affinoid subdomain isomorphic to a disjoint union of closed disks
  (set $C:=Y_\eta$ and $D:=]V[_Y$).
\end{enumerate}
\end{rem}

Analogous to the semistable reduction theorem for smooth projective
curves we have the following result. 

\begin{prop} \label{ssprop}
  Let $\Xs$ be an open analytic curve.  
  \begin{enumerate}
  \item After enlarging the field of definition, if necessary, there
    exists a semistable formal model of $\Xs$.
  \item If no connected component of $\Xs$ is a disk or an annulus
    then $\Xs$ has a stable formal model.
  \end{enumerate}
\end{prop}

\proof It is no restriction to assume that $\Xs$ is connected. By
definition we can write $\Xs=C-D$, where $C$ is a smooth projective
curve over $K$ and $D$ is a nonempty affinoid subdomain. After
enlarging $K$ we may assume that $C$ has a semistable model $Y$ over
$\OO$ (\cite{DeligneMumford69}). After an admissible blowup of $Y$
(which preserves semistability) we may further assume that the
affinoid $D\subset C$ is equal to the formal fiber $]V[_Y$ of an open
subscheme $V\subset Y_s$ (\cite{BoschLuetke93}). Set $Z:=Y_s-V$. Then
$\Xs= ]Z[_Y$ and therefore $Z$ is connected. In case $Z$ consists of a
single closed point, we do a further blowup to make sure that $Z$ is
the union of irreducible components of $Y_s$. Now let $\X$ be the
formal completion of $Y$ along $Z$. By construction $\X$ is a
semistable formal model of $\Xs$.  This proves (i).

Let $(Y',Z')$ denote the algebraic model of $\Xs$ obtained from
blowing down all instable components of $Z$. It is well known that the
curve $Y'$ is still semistable. There are two possible cases to
consider. The first case is that $Z'$ consists of a single point $z'$.
Then $z'$ is either a smooth point or an ordinary double point of
$Y_s'$. Therefore, $\Xs$ is either a disk or an annulus. In the second
case, $Z'$ is a union of stable irreducible components of $Y_s'$.
Therefore, $(Y',Z')$ is a stable algebraic model of $\Xs$ and the
formal completion of $Y'$ along $Z'$ is a stable formal model. This
finishes the proof of (ii).
\Endproof

\begin{prop} \label{minimalprop}
  Let $\Xs$ be an open analytic curve.
  \begin{enumerate}
  \item After a finite extension of the field of definition, there
    exists a minimal formal model of $\Xs$. It is unique up to unique
    isomorphism.
  \item The stable model of $\Xs$ (if it exists) is unique up to
    unique isomorphism.
  \end{enumerate}
\end{prop}

\proof By Proposition \ref{ssprop} (i) there exists a semistable
algebraic model $(Y,Z)$ of $\Xs$. Let $(Y',Z')$ denote the algebraic
model obtained from blowing down all irreducible components of $Z$ of
dimension one. Then $(Y',Z')$ is a minimal algebraic model. This
settles the existence part of (i). Let $U\subset Y'$ be an open affine
subset containing $Z'$. Let $\U$ denote the formal completion of $U$
along its special fiber. Since $\U$ is normal, it is the canonical
integral model of its generic fiber $\U_\eta$. Moreover,
$\Xs\subset\U_\eta$ is the disjoint union of the formal fibers of the
points contained in $Z'$. In this situation, a result of Bosch
\cite{Bosch77} (Korollar 5.9 and Satz 6.1) implies that $\X'$, the
formal completion of $\U$ in $Z'$ can be canonically identified with
$\Spf A$ where
\[
      A = H^0(\Xs,\OO_{\Xs}^\circ)
\]
is the ring of power bounded analytic functions on $\Xs$. This ring is
obviously independent of the choices we made and depends
functorially on $\Xs$. We conclude that the minimal formal model $\X'$
is unique up to unique isomorphism. This finishes the proof of (i).

To prove (ii), let $\X_1$ and $\X_2$ be stable models of $\Xs$. We
have to show that there exists a unique isomorphism $\X_1\cong\X_2$ of
formal schemes over $\OO$ extending the identity on $\Xs$. The
uniqueness of such an isomorphism is obvious.

For $i=1,2$ choose a stable algebraization $(Y_i,Z_i)$ of the formal
model $\X_i$. Using standard techniques, one shows that we may choose
the curve $Y_i$ to be stable. (This is not automatic, since Definition
\ref{modeldef} requires stability only for the irreducible components
contained in $Z_i$.) Let $(Y_i',Z_i')$ be the minimal algebraic model
obtained by blowing down all irreducible components of $Z_i$ (as in
the proof of (i)) and let $\X_i'$ denote the resulting minimal formal
model. By (i) we may identify $\X_1'$ with $\X_2'$. This means that we
can also identify the stable model $\X_2$ with the blowup of $\X_1'$
at an admissible sheaf of ideals $\hat{\J}$ on $\X_1'$ with support in
the point $z_1$. Let $\J$ denote the unique sheaf of ideals on $Y_1'$
with support in $z_1$ which gives rise to $\hat{\J}$ after completion
at $z_1$. Let $Y_3$ denote the blowup of $Y_1'$ in $\J$ and
$Z_3\subset Y_3$ the exceptional divisor.  By construction, the formal
completion of $Y_3$ along $Z_3$ can be identified with the formal
model $\X_2$.  It is also clear that the curve $Y_3$ is stable. Now
the uniqueness of the stable model for projective curves shows that
there exists a unique isomorphism $Y_1\cong Y_3$ extending the
identity on the generic fiber. It induces the desired isomorphism
$\X_1\cong\X_2$ of formal models.  \Endproof

\begin{rem}
  The proof of Proposition \ref{minimalprop} shows more generally that
  the minimal model of $\Xs$ is `minimal among all good models'.
\end{rem}

\subsection{The ends} \label{ends}

Let $\Xs$ be an open analytic curve.
The following definitions are taken over word by word from \cite{Coleman03}. 

\begin{defn} \label{opendef2} An {\em underlying affinoid} is an
  affinoid subdomain $U\subset\Xs$ such that $\Xs-U$ is the disjoint
  union of annuli none of which is contained in an affinoid subdomain
  of $\Xs$. An {\em end} of $\Xs$ is an element of the inverse limit
  of the set of connected components of $\Xs-U$, where $U$ ranges over
  all underlying affinoids. The set of all ends is denoted by
  $\partial\Xs$.
\end{defn}

By Remark \ref{formalrem} (v) we can write $\Xs=C-D$ where $C$ is a
smooth projective curve and $D=\cup_i D_i$ is a disjoint union of
closed disks $D_i$. If $U\subset \Xs$ is an underlying affinoid then
$C-U$ is disjoint union of open disks $D_i'$ such that $D_i\subset
D_i'$ and $D_i'-D_i$ is an open annulus representing an end of $\Xs$.
It follows that there is a natural bijection between the set of disks
$D_i$ and the set of ends $\partial\Xs$. See \cite{Coleman03}.

Let $(Y,Z)$ be a semistable algebraic model of $\Xs$ and $\X$ the
resulting formal model. We let $\partial Z\subset Z$ denote the subset
where $Z$ intersects an irreducible component of $Y_s$ not belonging
to $Z$ and call it the {\em boundary} of $Z$. The open subset
$Z^\circ:=Z-\partial Z$ is called the {\em interior} of $Z$. It is
clear that the definition of $\partial Z$ and $Z^\circ$ depend only on
the formal model $\X$ but not on the chosen algebraization.

Every $z\in \partial Z$ is a smooth point of $Z$. The formal fiber
$A_z:=]z[_{\X}\subset\Xs$ is an open annulus. Furthermore,
$U:=]Z^\circ[_{\X}\subset\Xs$ is an underlying affinoid for $\Xs$.
Hence we obtain a bijection between $\partial Z$ and $\partial\Xs$.
For every point $z\in\partial Z$ there is a unique closed formal
subscheme of the special fiber $T_z\subset\X_s$ with support in $z$
and isomorphic to $\Spf\, k[[t]]$. We call $T_z$ a {\em virtual
  component} of the reduction $Z$.  If $\xi$ denotes the end of $\Xs$
corresponding to $z$, we also write $T_\xi$ instead of $T_z$.

Let $\X\to\X'=\Spf A$ be the natural map onto the minimal formal model
which blows down $Z$ to a finite set of closed points (see the proof
of Proposition \ref{minimalprop}). This map identifies the virtual
components $T_\xi$ with the normalizations of the irreducible
components of the scheme $\Spec (A\otimes_{\OO} k)$. Therefore,
Proposition \ref{minimalprop} (i) implies that the scheme $T_\xi$
depends only, and in a functorial way, on the open analytic curve
$\Xs$ and the end $\xi$. In particular, any automorphism of $\Xs$
which fixes the end $\xi$ induces an automorphism of $T_\xi$. We say
that such an automorphism {\em acts trivially} on the end $\xi$ if it
induces the identity on $T_\xi$.

\section{Etale cohomology of open analytic curves}

\subsection{} \label{criterion1}

Let us fix, once and for all, a prime number $\ell$ which is prime to
$p$, the characteristic of the residue field $k$ of $K_0$. We also fix
an algebraic closure $\QQb_\ell$ of $\QQ_\ell$. Given an open analytic
curve $\Xs$ with field of definition $K$, we define
\[
    H^i(\Xs) := \bigl(\,\varprojlim_r H^i(\Xs\otimes_K\Kh_0\ac,\,\ZZ/\ell^r)\,\bigr)
       \otimes_{\ZZ_\ell}\QQb_\ell.
\]
Here $H^i(\,\cdot\,,\ZZ/\ell^r)$ is the \'etale cohomology of rigid
analytic spaces defined by Berkovich \cite{Berkovich93}. Similarly, we
define cohomology with compact support $H^i_c(\Xs)$ and, for every
affinoid subdomain $U\subset \Xs$, cohomology with support in $U$,
$H^i_U(\Xs)$. We define the {\em cuspidal part} of
cohomology as
\[
   H^i(\Xs)\cusp \,:=\, \text{image of \,}
              H^i_c(\Xs) \to H^i(\Xs).
\]
It follows from general facts that all these cohomology groups are
finite dimensional vector spaces over $\QQb_\ell$. We also have a
Poincar\'e duality isomorphism $H^i(\Xs)\cong H^{2-i}_c(\Xs)$.

By Remark \ref{formalrem} (v) we can write $\Xs=C-D$, where $C$ is a
smooth projective curve and $D=\cup_\xi D_\xi$ is a disjoint union of
closed disks corresponding to the ends of $\Xs$. This representation
induces a long exact cohomology sequence
\begin{equation} \label{longexactseq1}
    \cdots\to H^i_c(\Xs) \to H^i(C) \to \oplus_\xi H^i(D_\xi) \to \cdots
\end{equation}
Since $H^0_c(X)=H^1(D_\xi)=0$ we obtain a short exact sequence
\begin{equation} \label{exactseq1}
  0 \to B(\Xs) \to H^1_c(\Xs) \to H^1(C) \to 0,
\end{equation}
where the {\em boundary module} $B(\Xs)$ is defined by the short exact sequence\
\[
     0\to \QQb_\ell^{\pi_0(\Xs)} \to \QQb_\ell^{\partial\Xs} \to B(\Xs) 
      \to 0.
\]
In particular, if $\Xs$ is connected then $C$ is connected as well and we get
\[
       \dim H^1_c(\Xs) = 2g(C) +\abs{\partial\Xs} - 1.
\]
Taking the Poincar\'e dual of \eqref{longexactseq1} we get the long exact sequence
\[
    \cdots\to \oplus_\xi H^i_{D_\xi}(C) \to H^i(C) \to H^i(\Xs) \to \cdots
\]
which induces a short exact sequence 
\begin{equation} \label{exactseq2}
   0 \to H^1(C) \to H^1(\Xs) \to B(X)^* \to 0.
\end{equation}
Here $B(X)^*$ is the dual of $B(X)$. In particular, we obtain an isomorphism
\[
       H^1(\Xs)\cusp \cong H^1(C).
\]

\subsection{} \label{criterion2}

Let $\Xs$ be an open analytic curve, and let $\X$ be a good formal
model of $\Xs$. Let $Z$ denote the reduction of $\Xs$ with respect to
$\X$. Let $f:\X'\to\X$ be an admissible blowup such that $\X'$ is a
semistable model of $\Xs$ (exists after enlarging the field of
definition). Let $Z'$ be the reduction of $\X'$. This is a semistable
curve over $k$. At each of the finitely many points $z\in Z$ where $f$ is
not an isomorphism the inverse image $W_z:=f^{-1}(z)$ is a connected
union of irreducible components of $Z'$ and hence also a semistable
curve over $k$.

\begin{defn} \label{almostssdef} 
\begin{enumerate}
\item
  The model $\X$ is called {\em almost
  semistable} if the curves $W_z$ all have arithmetic genus zero.
\item The open analytic curve $\Xs$ is said to have {\em tree-like
    reduction} if the graph of components of the semistable curve $Z'$ is a tree.
\end{enumerate}
\end{defn}

By \cite{Berkovich96} the model $\X$ gives rise to
a (derived) sheaf of vanishing cycles on $Z$, which we denote by
$\RR\psi_{\X}$. The cohomology sheaves of $\RR\psi_{\X}$ are denoted
by $R^i\psi_{\X}$. By construction we have
\[
     (R^i\psi_{\X})_z = H^i(\Xs_z),
\]
where $\Xs_z:=]z[_{\X}\subset\Xs$ is the formal fiber of a
closed point $z\in Z$. (Since $\Xs_z$ is again an open analytic curve,
we see that the stalks of $R^i\psi_{\X}$ are finite dimensional and
that the passage from torsion coefficients to $\QQb_\ell$-coefficients
is justified.) Since for any closed point $z\in Z$ the formal fiber
$\Xs_z$ is connected (see Remark \ref{formalrem} (iv)) we have
$R^0\psi_{\X}=\QQb_\ell$. Therefore, the spectral sequence
$H^i(Z,R^j\psi_{\X})\Rightarrow H^{i+j}(\Xs)$ gives rise to the
exact sequence
\begin{equation} \label{vceq2}
   0 \to H^1(Z) \to H^1(\Xs) \to
    H^0(Z,R^1\psi_{\X}) \pfeil{d} H^2(Z).
\end{equation}

\begin{prop} \label{vcprop}
  The following two conditions are equivalent.
  \begin{enumerate}
  \item
    The model $\X$ is almost semistable and $\Xs$ has tree-like 
    reduction. 
  \item
    The first map in the sequence \eqref{vceq2} induces an isomorphism
    \[
          H^1(Z) \cong H^1(\Xs)\cusp.
    \]
  \end{enumerate}
\end{prop}

\proof Let $(Y',Z')$ be an algebraization of $\X'$ such that $Y'_s-Z'$
is a disjoint union of affine lines (see Remark \ref{formalrem} (v)).
Then $\Xs\subset Y_\eta'$ is the complement of closed disks and hence
we have a natural isomorphism $H^1(\Xs)\cusp\cong H^1(Y_\eta')$ (see
the previous subsection). Since $Y'$ is semistable, the
cospecialization map $H^1(Y_s)\to H^1(Y_\eta')$ is an isomorphism if
and only if the graph of components of $Y_s'$ is a tree. We conclude
that the natural map $H^1(Z')\to H^1(\Xs)$ induces an isomorphism
$H^1(Z')\cong H^1(\Xs)\cusp$ if and only if $\Xs$ has tree-like
reduction.

By \cite{Berkovich96}, Corollary 2.3 (ii), we have a natural
isomorphism of derived sheaves
\begin{equation} \label{fRpsi}
     \RR\psi_{\X} \cong f_*\RR\psi_{\X'}.
\end{equation}
In particular, we obtain a map between two exact sequences:
\begin{equation} \label{vcdiag}
\begin{array}{ccccccccc} 
  0 &\to& H^1(Z) &\to& H^1(\Xs) &\to& H^0(Z,R^1\psi_{\X}) &\pfeil{d}& H^2(Z) \\
    && \downarrow && \downarrow && \downarrow && \downarrow \\
  0 &\to& H^1(Z') &\to& H^1(\Xs) &\to& H^0(Z',R^1\psi_{\X'}) &\pfeil{d}& H^2(Z'). \\
\end{array}
\end{equation}
In this diagram, the first, the second and the fourth vertical arrows
are injective.  The third vertical arrow is in general not injective;
the isomorphism \eqref{fRpsi} induces an isomorphism 
\begin{equation} \label{kerWz}
   \Ker\big(H^0(Z,R^1\psi_{\X})\to H^0(Z',R^1\psi_{\X'})\big) \cong
        \oplus_z \, H^1(W_z),
\end{equation}
where $z\in Z$ runs over the points where $f$ is not an isomorphism.
A simple diagram chase yields an isomorphism of \eqref{kerWz} with
the cokernel of the first vertical map in \eqref{vcdiag}.

Now suppose that (ii) holds. Then the discussion in the first
paragraph of this proof, together with the injectivity of $H^1(Z)\to
H^1(Z')$, shows that $\Xs$ has tree-like reduction.  Furthermore,
$H^1(Z)\cong H^1(Z')$ and hence $H^1(W_z)=0$ for all critical points
$z$. Since $W_z$ is semistable, this means that all $W_z$ have genus
zero. We have proved that (ii) implies (i). The proof of the converse
is similar.  \Endproof

\begin{rem}
  Let $\X$ be a good formal model of $\Xs$ whose reduction $Z$ is
  purely of dimension one. Then $\X$ is semistable if and only if for
  every closed point $z\in Z$ we have
  \[
            \dim (R^1\psi_{\X})_z \leq 1.
  \]
  Indeed, the above condition implies that the formal fiber $\Xs_z$ is
  either a disk or an annulus.
\end{rem}

\section{Etale Galois covers} \label{diskcover}

In this section we fix an open analytic curve $\Ys$, a finite group
$G$ and an \'etale $G$-torsor $f:\Xs\to\Ys$. Let $K$ denote the field
of definition of $\Ys$ (assumed to be sufficiently large). We assume
that the auxiliary prime $\ell$ chosen in the last section does not
divide the order of $G$.

\subsection{Algebraization} \label{diskcover1}

By Remark \ref{formalrem} (v) we can write $\Ys=C-D$, where $C$ is a
smooth projective curve over $K$ and $D$ is an affinoid subdomain,
isomorphic to a finite union of closed disks.

\begin{prop} \label{Marcoprop}
\begin{enumerate}
\item After a finite extension of $K$, the \'etale
  $G$-torsor $f$ extends to a finite (possibly ramified) $G$-cover $C'\to
  C$ of smooth projective curves over $K$.
\item
  The rigid space $\Xs$ is an open analytic curve.
\end{enumerate}
\end{prop}

\proof Part (i) follows from a result of Garuti \cite{Garuti96}. Part
(ii) is an immediate consequence of (i), because the inverse image of
the affinoid $D$ in $C'$ is again an affinoid.  \Endproof

\subsection{Formal models} \label{diskcover2}

A {\em semistable} (resp.\ a {\em minimal}) {\em formal model} of the
$G$-torsor $f:\Xs\to\Ys$ is given by a finite, $G$-invariant morphism
of formal schemes $\X\to\Y$ extending $f$, where $\X$ and $\Y$ are
semistable (resp.\ minimal) models of $\Xs$ and $\Ys$.

We claim that there always exists a minimal and a semistable model of
$f$.  For the minimal model this is easy: the minimal model $\X$ of
$\Xs$ exists and is unique (Proposition \ref{minimalprop} (i)).
Therefore, the $G$-action extends to $\X$ and we can take $\Y:=\X/G$.
To obtain the semistable model, let us first assume that $\Xs$ has a
stable model $\X$. Again we can use unicity to extend the $G$-action
and set $\Y:=\X/G$. That $\Y$ is semistable follows from
\cite{Raynaud90}, Proposition 5. If $\Xs$ does not have a stable model
then its connected components are open disks or annuli (Proposition
\ref{ssprop}). The claim in this case is left as an exercise.

Conversely, let us start with a good formal model $\Y$ of $\Ys$. We
claim that there exists a unique formal model $\X$ of $\Xs$ which is
normal and such that the $G$-torsor $f:\Xs\to\Ys$ extends to a finite
morphism $\X\to\Y$ and we have $\Y=\X/G$.  Indeed, let $(Y,Z)$ be an
algebraization of the formal model $\Y$ as in Remark \ref{formalrem}
(v).  By Proposition \ref{Marcoprop} (ii), the $G$-torsor $f$ extends
to a finite cover $C'\to C=Y_\eta$ of smooth projective curves. Let
$Y'$ be the normalization of $Y$ in $C'$ and let $Z'\subset Y_s'$
denote the inverse image of $Z\subset Y_s$. Then $(Y',Z')$ is an
algebraic model of $\Xs$ and gives rise to the desired formal model
$\X$. To show that $\X$ is unique use \cite{Fargues05}, Appendix A.

In general, the
formal model $\X$ will not be good because its special fiber may not
be reduced. However, after replacing the field of definition $K$ by a
finite extension, we may assume that $\X_s$ is reduced and so $\X$ is
good. This follows easily from the Reduced Fiber Theorem of Grauert
and Remmert (\cite{BGR}, \S 6.4.1). By Remark \ref{formalrem} (iii) the
formation of $\X$ is then stable under further extension of $K$. Hence
the definition of the normalization is compatible with our philosophy
of keeping the field of definition $K$ variable and sufficiently
large.

\subsection{Decomposition and inertia groups} \label{diskcover3}

Let $U\subset \Ys$ be an affinoid subdomain. The canonical reduction
of $U$ is an affine curve over $k$ which we denote by $\bar{U}$. By
the Reduced Fiber Theorem \cite{BGR} we may and will assume that
$\bar{U}$ is reduced (after replacing the field of definition by a
suitable finite extension). We say that $U$ has {\em good} (resp.\
{\em irreducible}) {\em reduction} if $\bar{U}$ is smooth over $k$
(resp.\ irreducible).

Let $U\subset\Ys$ be an affinoid with good and irreducible reduction.
Choose a connected component $V$ of the inverse image $f^{-1}(U)$.
Clearly, $V$ is an affinoid subdomain of $\Xs$. We say that $U$ has
{\em good} (resp.\ {\em irreducible}) {\em reduction} in $\Xs$ if $V$
has good (resp.\ irreducible) reduction. We note that the canonical
reduction $\bar{V}$ is connected because $V$ is connected. Therefore,
if $U$ has good reduction in $\Xs$ then it also has irreducible
reduction in $\Xs$.

Suppose that $U$ has irreducible reduction in $\Xs$. The stabilizer in
$G$ of the connected component $V$ is denoted by $G(U)$ and is called
the {\em decomposition group} of $U$ (it is independent of the choice
of $V$, up to conjugation in $G$). The {\em inertia group} $I(U)\lhd
G(U)$ of $U$ is defined as the kernel of the natural homomorphism
\[
       G(U) \To\Aut_k(\bar{V}).
\]

\begin{prop} \label{irredredprop}
  Let $U\subset\Ys$ be an affinoid with good and irreducible
  reduction. Let $f:\Xs\to\Ys$ be an \'etale $G$-torsor.
  \begin{enumerate}
  \item Suppose that $U$ has irreducible reduction in $\Xs$. Choose a
    connected component $V\subset f^{-1}(U)$ with stabilizer $G(U)$
    and set $V':=V/I(U)$. Then the natural map 
    \[
            \bar{V}'\to\bar{U}
    \]
    is an \'etale Galois cover of smooth affine curves over $k$, with
    Galois group $G(U)/I(U)$. Furthermore, the natural map
    \[
           \bar{V}\to\bar{V}'
    \]
    is finite and radicial of degree $|I(U)|$. Since $\bar{V}$ is
    reduced it follows that $I(U)$ is a $p$-group (where $p$ is the
    characteristic of $k$).
  \item Suppose that $G$ is a $p$-group. Then $U$ has irreducible
    reduction in $\Xs$. 
  \end{enumerate}
\end{prop}

\proof Part (i) is a consequence of the Purity Theorem of
Zariski-Nagata and the assumption that $f$ is \'etale. See e.g.\
\cite{Raynaud99}, \S 2.4. Under the assumption that $G$ is a $p$-group
it is proved in \cite{Raynaud90} that the singularities of $\bar{V}$
are unibranched. Therefore, connectedness of $\bar{V}$ implies its
irreducibility.  \Endproof

Let $\xi\in\partial\Ys$ be an end of $\Ys$. Choose an end $\xi'$ of
$\Xs$ lying above $\xi$ and let $G(\xi)\subset G$ denote the
stabilizer of $\xi'$ (which depends, up to conjugation in $G$, only on
$\xi$). We call $G(\xi)$ the {\em decomposition group} of $\xi$ in
$\Xs$. The {\em inertia group} $I(\xi)\lhd G(\xi)$ of $\xi$ is the
subgroup of elements which act trivially on $\xi$ (see \S \ref{ends}).

\subsection{Vanishing cycle sheaves for nonconstant coefficients} \label{diskcover4}

Let $\tau:G\to\GL(W)$ be a representation of $G$ on a finite
dimensional $\QQb_\ell$-vector space $W$. Given a subgroup $H\subset
G$ which acts on a $\QQb_\ell$-vector space $W'$, we set
$W'[\tau]:=\Hom_{H}(W,W')$ (which subgroup $H$ we mean should always
be clear from the context). Extending this definition from
$\QQb_\ell$-vector spaces to sheaves, we obtain an exact functor
$\F\mapsto (f_*\F)[\tau]$ from $\QQb_\ell$-sheaves on $\Xs$ with
$G$-action to $\QQb_{\ell}$-sheaves on $\Ys$. (To construct this
functor, one has to choose a finite extension $E/\QQb_\ell$ and a
projective $\OO_E[G]$-module $M$ such that $W=M\otimes\QQb_\ell$, and
define everything first for $\OO_E/\ell^n$-sheaves.  Since we assume
that $\ell$ does not divide the order of $G$, this poses no problem.)
As a special case of this construction we set
\begin{equation} \label{Fdefeq}
     \F_\tau := (f_*\QQb_\ell)[\tau].
\end{equation}
General arguments show that
\[
   H^i(\Xs,\QQb_\ell)[\tau] = H^i(\Ys,\F_\tau).
\]
A similar equality holds for cohomology with support and for the
cuspidal part $H^i(\Ys,\F_\tau)\cusp$ (which is defined as the image of the
map $H^i_c(\Ys,\F_\tau)\to H^i(\Ys,\F_\tau)$). 

Let $U\subset\Ys$ be an affinoid subdomain. By \cite{Berkovich96} $\F_\tau$
gives rise to a derived sheaf of vanishing cycles
$\RR\psi\F_\tau|_{\bar{U}}$ on the canonical reduction $\bar{U}$ such that
\begin{equation}
   H^i_U(\Ys,\F_\tau) \cong \mathbb{H}^i_c(\bar{U},\RR\psi\F_\tau|_{\bar{U}}).
\end{equation}
In particular, we obtain an exact sequence
\begin{multline} \label{Fexactseq}
  \qquad 0 \to H^1_c(\bar{U},R^0\psi\F_\tau|_{\bar{U}}) \to H^1_U(\Ys,\F_\tau) \to
        H^0_c(\bar{U},R^1\psi\F_\tau|_{\bar{U}}) \\ \to 
          H^2_c(\bar{U},R^0\psi\F_\tau|_{\bar{U}}). \qquad
\end{multline}
Let $V\subset\Xs$ be a connected component of $f^{-1}(U)$. Let
$\varphi:\bar{V}\to\bar{U}$ denote the finite morphism induced by the
natural map $V\to U$. Using \cite{Berkovich96}, Corollary 2.3 (ii) and
the fact that $\varphi_*$ is exact, we obtain a canonical isomorphism
\begin{equation} \label{equiiso0}
       \RR\psi\F_\tau|_{\bar{U}} \cong \varphi_*(\RR\psi\QQb_\ell|_{\bar{V}})[\tau].
\end{equation}
In particular, we have a canonical isomorphism
\begin{equation} \label{equiiso}
    R^0\psi\F_\tau|_{\bar{U}} \cong (\varphi_*\QQb_\ell)[\tau].
\end{equation}

\begin{defn} \label{resolvedef}
  Let $U\subset\Ys$ be an affinoid subdomain. We say that the sheaf
  $\F_\tau$ is {\em residual} over $U$ if the first map in
  \eqref{Fexactseq} is an isomorphism,
  \[
        H^1_c(\bar{U},R^0\psi\F_\tau|_{\bar{U}}) \iso H^1_U(\Ys,\F_\tau).
  \]
  We say that $\F_\tau$ is {\em resolved} over $U$ if it is
  residual over $U$ and the natural map 
  \[
         H^1_U(\Ys,\F_\tau)\to H^1(\Ys,\F_\tau)
  \]
  is surjective.
\end{defn}

It is of course a nontrivial problem find an affinoid over which the
sheaf $\F_\tau$ is resolved. We limit the discussion of this
problem to the following proposition which gives a criterion for
$\F_\tau$ to be residual, and we refer the reader to \cite{LT} for
concrete and nontrivial examples.

\begin{prop}
  The sheaf $\F_\tau$ is residual over $U$ if either one of the following conditions holds.
  \begin{enumerate}
  \item
    The affinoid $U$ has good reduction in $\Xs$. 
  \item
    The affinoid $U$ has irreducible reduction, and the
    restriction of $\tau$ to the inertia group $I(U)\subset G$ is
    trivial. 
  \end{enumerate}
\end{prop}

\proof If $U$ has good reduction in $\Xs$ then the vanishing cycle
sheaves $R^i\psi\F_\tau|_{\bar{U}}$ are zero for $i>1$, by \eqref{equiiso0}.
Therefore, $\F_\tau$ is residual over $U$. 

Now assume that Condition (ii) holds.  Let $V\subset\Xs$ be a
connected component of $f^{-1}(U)$ with decomposition group
$G(U)\subset G$. It follows from \cite{AbbesSaito02}, \S 5, that there
exists an open analytic curve $Y_1\subset Y$ containing $U$ such that
$f^{-1}(Y_1)$ has a {\em unique} connected component $\Xs_1\subset\Xs$
containing $V$. Set $\Xs':=\Xs_1/I(U)$ and $V':=V/I(U)$. Clearly,
$f':\Xs'\to\Ys_1$ is an \'etale $G(U)/I(U)$-torsor.  Our assumption
says that the restriction of $\tau$ to $G(U)$ comes from a
representation $\tau'$ of $G(U)/I(U)$. Hence the restriction of
$\F_\tau$ to $\Ys_1$ is isomorphic to the sheaf
$\F_{\tau'}:=(f'_*\QQb_\ell)[\tau']$. By excision we obtain an
isomorphism
\[
    H^1_U(\Ys,\F_\tau) \cong H^1_U(\Ys_1,\F_{\tau'}).
\]
Furthermore, it follows from Proposition \ref{irredredprop} (i) that
$U$ has good reduction in $\Xs'$. Therefore, we can use the first case
of the proposition which is already proved.
\Endproof

\subsection{An equivariant criterion for almost semistability} \label{diskcover5}

Let $\Y$ be a semi\-stable formal model of $\Ys$. We define the formal
model $\X$ of $\Xs$ as the normalization of $\Y$ in $\Xs$ (see \S
\ref{diskcover2}). We would like to have a criterion that ensures that
$\X$ is semistable. We will denote the reduction of $\X$ (resp.\ of
$\Y$) by $Z$ (resp.\ $Z'$).

\begin{defn} \label{suppdef} Let $U\subset\Ys$ be an affinoid
  subdomain. We say that $U$ is {\em supported} by the semistable
  model $\Y$ if there exists an open subset $W\subset Z'$
  such that $U=]W[_{\Y}$.
\end{defn}

Suppose that $U$ is supported by $\Y$, and let $W\subset Z'$ be as in
the definition. Let $\U$ be the canonical integral model of $U$ (with
special fiber $\bar{U}$) and $\W\subset\Y$ the open formal subscheme
whose underlying topological space is $W$. Since $\W$ is a normal
model of the affinoid $U$, there exists a unique morphism of formal
schemes $\W\to\U$ extending the identity on the generic fiber $U$.  In
fact, the morphism $\W\to\U$ is an admissible blowup. On the special
fiber we obtain a proper and surjective morphism $W\to\bar{U}$ which
is an isomorphism over some dense open subset of $\bar{U}$. See
\cite{BoschLuetke93}. 

\begin{rem} \label{supprem} The following facts are easy
  consequences of Definition \ref{suppdef} and the discussion
  following it.
  \begin{enumerate}
  \item If $U$ has good reduction and is supported by $\Y$ then the
    morphism $W\to\bar{U}$ has a section. In particular, we have a
    canonical locally closed embedding $\bar{U}\inj Z'$.
  \item Suppose $U$ is supported by $\Y$. Let $V\subset\Xs$ be a
    connected component of $f^{-1}(U)$. Then the affinoid $V$ is
    supported by the model $\X$.
  \item Given a finite family of affinoids $U_i\subset\Ys$ which
    meets every connected component of $\Ys$, there exists a minimal
    semistable model $\Y$ of $\Ys$ which supports all $U_i$.
  \end{enumerate}
\end{rem}

\begin{prop} \label{criterionprop}
  Suppose that for every irreducible $G$-representation $\tau$ such that 
  \[
       H^1(\Ys,\F_\tau)\cusp \neq 0
  \]
  there exists an affinoid $U_\tau\subset\Ys$ over which $\F_\tau$ is
  resolved. Let $\Y$ be a semistable model of
  $\Ys$ which supports $U_\tau$ for every $\tau$. Let $\X$ be the
  normalization of $\Y$ in $\Xs$. Then $\X$ is almost semistable.
  Moreover, $\Xs$ has tree-like reduction.
\end{prop}

\proof Let $\Y$ and $\X$ be the formal models from the statement of
the proposition. Let $Z$ denote the reduction of $\X$. By Proposition
\ref{vcprop} it suffices to show that we have an isomorphism
\[
      H^1(Z) \cong H^1(\Xs)\cusp.
\]
Actually, since the natural map $H^1(Z)\to H^1(\Xs)$ is injective, it
suffices to show that $H^1(Z)$ maps onto $H^1(\Xs)\cusp$. The
group $G$ acts on both these vector spaces, hence it suffices to show that
for every irreducible $G$-representation $\tau$ the map $H^1(Z)\to
H^1(\Xs)$ induces a surjective map
\begin{equation} \label{criterioneq1}
    H^1(Z)[\tau] \to H^1(\Xs)\cusp[\tau]= H^1(\Ys,\F_\tau)\cusp.
\end{equation}
We may assume that the right hand side of \eqref{criterioneq1} is not
zero. Therefore, our hypothesis says that there exists an affinoid
$U\subset\Ys$ which is supported by $\Y$ and over which the sheaf
$\F_\tau$ is resolved. Let $V\subset\Xs$ be a connected component of
$f^{-1}(U)$. By Remark \ref{supprem} (ii) there exists an open subset
$W\subset Z$ and a proper surjective map $W\to\bar{V}$. We thus obtain
natural maps
\[
       H^1_c(\bar{V}) \to H^1_c(W) \to H^1(Z).
\]
By \eqref{equiiso}, the composition of these two maps induces the left
vertical arrow in the following diagram.
\begin{equation} \label{criterioneq2}
\begin{CD}
  H^1_c(\bar{U},R^0\psi\F_\tau|_{\bar{U}}) @>{\cong}>> H^1_U(\Ys,\F_\tau)    \\
  @VVV                       @VVV                 \\
  H^1(Z)[\tau]         @>>> H^1(\Ys,\F_\tau)\cusp \\
\end{CD}
\end{equation}
The upper horizontal arrow is an isomorphism because $\F_\tau$ is
residual over $U$.  The right vertical arrow is surjective because
$\F_\tau$ is resolved over $U$.  It follows that the lower
horizontal arrow is also surjective. This finishes the proof of the
proposition.  \Endproof

\begin{rem} 
  The criterion given by Proposition \ref{criterionprop} is sharp, in
  the following sense. Suppose that $\Xs$ has tree-like reduction.
  Then there exists an affinoid $U\subset\Ys$ such that $\F_\tau$ is
  resolved over $U$, for all $G$-representations $\tau$.

  To construct $U$, let $\X\to\Y$ be a semistable model of
  $f:\Xs\to\Ys$ and let $Z$ denote the reduction of $\X$. The complement
  \[
        W:=Z-\partial Z - Z^{\rm sing}
  \]
  is an affine and dense open subset. One checks that the natural map
  $H^1_c(W)\to H^1(Z)$ is surjective. Set $V:=]W[_{\X}$ and set
  $U:=f(V)$. By construction, $U$ is an affinoid over which $\F_\tau$
  is resolved for every $G$-representation $\tau$.

  It is clear that in practice this construction is not very useful to
  determine a stable model of $f$. The point of Proposition
  \ref{criterionprop} is rather that it suffices to look for a
  suitable affinoid for each irreducible representation $\tau$ at a
  time. Also, such an affinoid may be much simpler, and neither include
  nor be included in the affinoid $U$ constructed above.
\end{rem}

\subsection{Etale covers of the disk} \label{diskcover6}

In this final subsection we assume that $\Ys$ is an open disk. In this
case we can describe the semistable model $\Y$ in terms of a {\em tree
  of disks}. By a {\em closed disk} we shall mean an affinoid subdomain of $\Ys$
which is isomorphic (possibly after enlarging the field of definition)
to the closed unit disk.

Let $S$ be a nonempty finite collection of closed
disks. For every nonempty subset $T\subset S$ there exists a closed
disk $D_{T}\subset\Ys$ which is minimal with the property that
$D\subset D_{T}$ for every $D\in T$.  We say that $S$ is {\em closed}
if $D_{T}\in S$ for every nonempty subset $T\subset S$. It is easy to
see that for every finite collection of closed disks $S$ there is a
minimal finite collection of closed disks $\bar{S}$ which contains $S$
and is closed. We call $\bar{S}$ the {\em closure} of $S$.

Suppose that $S$ is closed. We will associate to $S$ a directed graph
$\Gamma=\Gamma_{S}$, as follows. Write $S=\{D_v\mid v\in V\}$ for an
index set $V$, and consider the elements of $V$ as vertices of
$\Gamma$. We add to $\Gamma$ a distinguished vertex $v_0$, called the
{\em boundary}. For every $v\in V$, the disk $D_v$ is either maximal
among the disks in $S$, or there exists $w\in V$ such that the disk
$D_w\in S$ is minimal among all disks in $S$ strictly containing
$D_v$. In the first case, we add to $\Gamma$ the edge $(v_0,v)$. In
the second case, we add the edge $(w,v)$. It is clear that $\Gamma$ is
a rooted and connected tree and that the boundary $v_0$ is the root.
The datum $(\Gamma;D_v)$ is called a {\em tree of disks} in $\Ys$. One
easily shows:

\begin{prop}
  For every tree of disks $(\Gamma;D_v)$ there exists a semistable
  model $\Y$ of $\Ys$ which is minimal with the property that every
  disk $D_v$ is supported by $\Y$. Let $Z'$ denote the
  reduction of $\Y$. Then the graph of components of the semistable
  curve $Z'$ is naturally isomorphic to the tree $\Gamma$ (the
  boundary $v_0$ corresponds to the virtual component $T_\xi$, where
  $\xi$ is the unique end of $\Ys$).  Furthermore, every semistable
  model of $\Ys$ arises in this way from a unique tree of disks.
\end{prop}

Let us fix a tree of disks $(\Gamma;D_v)$, and let $\Y$ be the
corresponding semistable model. Let $f:\Xs\to\Ys$ be the \'etale
$G$-cover fixed at the beginning of this section. Let $\X$ be the
normalization of $\Y$ in $\Xs$. We obtain the following reformulation
of Proposition \ref{criterionprop}.

\begin{prop} \label{criterionprop2}
  Suppose that for every irreducible $G$-representation $\tau$ such that 
  \[
       H^1(\Ys,\F_\tau) \neq 0
  \]
  there exists a set $V'\subset V$ of vertices of $\Gamma$ such that
  $\F_\tau$ is resolved over the affinoid
  \[
         U := \bigcup_{v\in V'} D_v.  
  \]
  Then the formal model $\X$ is almost semistable. Moreover, $\Xs$ has
  tree-like reduction.
\end{prop}

\begin{rem}
  In contrast to its ancestor, Proposition \ref{criterionprop}, this
  criterion is not sharp at all. In fact, given an \'etale $G$-Galois
  cover of the disk with tree-like reduction, there is no reason to
  expect that one can always resolve the sheaf $\F_\tau$ over a union
  of closed disks. Nevertheless, this is true in important examples,
  for instance for the Lubin-Tate spaces studied in \cite{LT}.
\end{rem}

For the rest of this section we assume that the conclusion of
Proposition \ref{criterionprop2} holds. We wish to give a further
criterion which ensures that the model $\X$ is semistable (and not just almost semistable). 

We fix the following notation. Let $Z'$ denote the reduction of $\Y$
and $Z$ the reduction of $\X$.  Let $e=(v_1,v_2)$ be an edge of the
graph $\Gamma$. To $e$ corresponds a closed point $z'\in Z'$ such that
$A_e:=]z'[_{\Y}$ is an open annulus. If $v_1\neq v_0$ then $z'$ is an
ordinary double point of $Z'$ and $A_e=D_e-D_{v_2}$, where $D_e\subset
D_{v_1}$ is the formal fiber containing $D_{v_2}$ (an open disk). If
$v_1=v_0$ then $z'$ is the unique element of the boundary $\partial
Z'$ of $Z'$.  In this case we have $A_e=\Ys-D_{v_2}$ and we set
$D_e:=\Ys$. For $i=1,2$ let $\xi_i'$ be the end of $A_e$ corresponding
to the vertex $v_i$.

We choose a point $z\in Z$ lying over $z'$ and set $\Xs_z:=]z[_{\X}$.
Note that $\Xs_z$ is a connected component of the inverse image of the
annulus $A_e$.  Choose, for $i=1,2$, an end $\xi_i\in\partial\Xs_z$
lying over $\xi_i'$. Since the formal model $\X$ is almost semistable,
$\Xs_z$ is an open analytic curve of genus zero. We want to show that
$\Xs_z$ is actually an open annulus.  The last part of the following
proposition gives a criterion when this is true.

\begin{prop} \label{criterionprop3}
  Suppose that $G$ is a $p$-group. Then the following holds.
  \begin{enumerate}
  \item The closed disks $D_v\subset\Ys$ have irreducible reduction in
    $\Xs$. In particular, the decomposition groups $G(D_v)\subset G$ and
    inertia groups $I(D_v)\lhd G(D_v)$ are well defined.
  \item Up to conjugation in $G$ we have
    \[
        G(D_{v_2})=G(\xi_2)\subset I(\xi_1)=I(D_{v_1}).
    \]
  \item The open curve $\Xs_z$ is an open annulus if and only if 
    \begin{equation} \label{G=Ieq}
          G(D_{v_2})= I(D_{v_1}).
    \end{equation}
  \end{enumerate}
  In particular, if \eqref{G=Ieq} holds for all edges $e=(v_1,v_2)$ then the
  formal model $\X$ is semistable.
\end{prop}

\proof Part (i) is a special case of Proposition \ref{irredredprop}
(ii). Part (ii) follows from the following lemma, applied to the
Galois cover $C_e\to D_e$, where $D_e\subset D_{v_1}$ is the residue
class containing the disk $D_{v_2}$ if $v_1\neq v_0$ (and $D_e:=\Ys$
otherwise) and $C_e$ is a connected component of $f^{-1}(D_e)$. Part
(iii) is now clear.  \Endproof

\begin{lem} \label{disklem} Let $G$ be a $p$-group, $\Ys$ an open disk
  and $f:\Xs\to\Ys$ an \'etale $G$-torsor. Suppose that $\Xs$ is
  connected.  Then the
  following holds.
  \begin{enumerate}
  \item
     The open analytic curve $\Xs$ has a unique end.
   \item Let $D\subset \Ys$ be a closed disk. Set $A:=\Ys-D$ and let
     $\xi_1$, $\xi_2$ denote the two ends of the open annulus $A$. We
     assume that $\xi_1$ corresponds to the unique end of $\Ys$. Then
     we have (up to conjugation in $G$):
     \[
         G = I(\xi_1) \supset G(\xi_2) = G(D).
     \]
   \end{enumerate}
\end{lem}

\proof (compare with \cite{Raynaud90}) Let $D\subset \Ys$ be a closed
disk and $\Y$ the minimal semistable model of $\Ys$ supporting $D$.
Let $Z'$ denote the reduction of $\Y$ and $z'\in\partial Z'$ the
unique boundary point. Then $Z'\cong\PP^1_k$, and
$Z'-\{z'\}\cong\AA^1_k$ is the canonical reduction of the disk $D$.
Let $\X$ be the normalization of $\Y$ in $\Ys$ and $Z$ the reduction
of $\X$. Let $C\subset f^{-1}(D)$ be a connected component. Then the
canonical reduction $\bar{C}$ of the affinoid $C$ can be identified
with an open subset of $Z$. By Proposition \ref{irredredprop} the
curve $\bar{C}$ is irreducible.  Furthermore, the map
$\bar{C}\to\bar{D}=Z'-\{z'\}$ factors as the composition of a finite
radicial map $\bar{C}\to\bar{C}'$ and an \'etale Galois cover
$\bar{C}'\to\bar{D}\cong \AA^1_k$ whose Galois group is a $p$-group. A
well known lemma says that such a cover of the affine line is totally
ramified over infinity. (One uses the fact that any proper subgroup of
a $p$-group is contained in a proper normal subgroup.) We conclude
that there exists a unique point $z\in Z$ which is mapped to $z'\in
Z'$ and lies in the closure of $\bar{C}$.  Furthermore, there exists a
unique branch of $Z$ through the point $z$ whose generic point lies on
$\bar{C}$.  Since the ends of $f^{-1}(A)$ lying above $\xi_2$ are in
natural bijection with the branches of $Z$ through $z$, the equality
\begin{equation} \label{plemeq1}
       G(\xi_2) = G(D)
\end{equation}
follows immediately.

Let $\X'\to\Y'$ be a semistable model of $f:\Xs\to\Ys$ and apply the
previous argument to the largest disk $D$ which is supported by the
model $\Y'$. Then by construction, $f^{-1}(D)$ is an underlying
affinoid of $\Xs$. In particular, $f^{-1}(D)$ is connected and
$f^{-1}(A)$ is a disjoint union of open annuli. Now \eqref{plemeq1}
implies that $f^{-1}(A)$ is connected and that
\begin{equation} 
       G=G(\xi_1).
\end{equation}
Part (i) of the proposition follows. To finish the proof of (ii) we
have to show that $I(\xi_1)=G(\xi_1)$. After dividing out by
$I(\xi_1)$ (which is a normal subgroup of $G=G(\xi_1)$!) we may assume
that $I(\xi_1)=1$. Applying purity of branch locus to the minimal
formal model $\X''\to\Y''=\Spf\,\OO[[t]]$ of $f:\Xs\to\Ys$, we conclude that $G=1$.
This finishes the proof of the proposition.
\Endproof

\begin{cor}
  Let $G$ be a finite $p$-group and $f:\Xs\to\Ys=D(0,1)$ be an \'etale
  Galois cover of the open unit disk. Then 
  \[
        H^1(\Xs) = H^1_c(\Xs) = H^1(\Xs)\cusp.
  \]
\end{cor}



\vspace{1cm}
\noindent
Stefan Wewers\\[1ex]
Mathematisches Institut\\
Beringstr. 1\\
73115 Bonn\\[1ex]
wewers@math.uni-bonn.de

\end{document}